\documentclass[12pt]{amsart}
\usepackage{amsfonts}

\theoremstyle{definition}

\theoremstyle{remark}

\newcommand{\const}{\mathop{\rm const}\limits}

\begin{document}

\begin{center}

{\bf  HARDY-LITTLEWOOD INEQUALITIES FOR RIESZ'S POTENTIAL: \\
low bounds estimations for different powers.} \par

\vspace{3mm}
{\bf E. Ostrovsky}\\

e - mail: galo@list.ru \\

\vspace{3mm}

{\bf L. Sirota}\\

e - mail: sirota@zahav.net.il \\

\vspace{4mm}

 Abstract. \\
{\it In this article we obtain the non - asymptotical low estimations for bilinear
Riesz's functional through the Lebesgue spaces norms by means of building of
 some examples.} \\

\end{center}

\vspace{3mm}


\vspace{3mm}


\vspace{3mm}

\section{Introduction. Statement of problem.}

\vspace{3mm}

 The linear integral operator $ I_{\alpha} f(x), $ or, wore precisely, the
 {\it family of operators} of a view

 $$
u(x) = I_{\alpha}f(x) = \int_{R^d} \frac{f(y) \ dy}{|x - y|^{d - \alpha}}
 $$
is called {\it Riesz's integral operator,} or simply {\it Riesz's} potential, or
{\it fractional integral.} \par
 The bilinear {\it Riesz's functional} $ B_{\alpha,d}(f,g) = B_{\alpha}(f,g)
= B(f,g) $ may be defined as follows:

 $$
 B(f,g) = \int_{R^d} \int_{R^d} \frac{f(y) \ g(x) \ dx \ dy}{|x - y|^{d - \alpha}}.
 $$
 It is evident that

 $$
 B_{\alpha,d}(f,g) = (I_{\alpha}f, g),
 $$
where $ (f,g)$ denotes ordinary inner product of a (measurable) functions $ f $
and $ g: $

$$
(f,g) = \int_{R^d} f(x) \ g(x) \ dx,
$$
which is defined, e.g. when $ f \in L_p, \ g \in L_q, \ p,q > 1, 1/p + 1/q = 1,
L_p = L_p(R^d) $ is the classical Lebesgue space of all the measurable functions
$ f: R^d \to R $ with finite norm

 $$
 |f|_p \stackrel{def}{=} \left( \int_{R^d} |f(x)|^p \ dx  \right)^{1/p};  \  f \in L_p \
 \Leftrightarrow |f|_p < \infty.
 $$
  Obviously,

  $$
  \sup_{|g|_q =1}  |B_{\alpha,d}(f,g)| = |I_{\alpha}f|_p.
  $$

Here  $ \alpha = \const \in (0,d) $  and $  |y| $ denotes usually  Euclidean norm
of the finite dimensional vector $ y. $  In particular, if $ y \in R^1,  \ |y|$
denotes ordinary absolute value of the number $ y. $ \par
In the case $ d \ge 2 $ and $ \alpha = 2 \ u(x) $ coincides with the classical
Newton's potential.\par
 The operators  $ I_{\alpha} $ and functionals $ B_{\alpha} $
 are used in the theory of Fourier transform, theory of Partial Differential
 Equations, probability theory (study of potential functions for Markovian processes
  and spectral densities for stationary random fields), in the
 functional analysis, in particular, in the theory of interpolation of operators etc.,
 see for instance \cite{Bennet1}, \cite{Stein1}, \cite{Lieb0},\cite{Lieb2}, \cite{Rubin1},
 \cite{Opic1} etc. \par

 We  denote also $ L(a,b) = \cap_{p \in (a,b)} L_p.  $ \par

{\bf  We will investigate the estimations of a view:}

$$
|B_{\alpha,d}(f,g) \le K_{\alpha,d}(r,s) \ |f|_r \ |g|_s. \eqno(1)
$$
 It is known \cite{Lieb1}, \cite{Lieb3}, \cite{Talenti1}
 that the inequality (1) is possible only in the case when
 $$
 r,s > 1, \ \frac{1}{r} + \frac{1}{s} = 1 + \frac{\alpha}{d}. \eqno(2)
 $$
We will denote the set of all such a values $ r,s, r,s > 1 $ as $ G = G_{\alpha,d}. $ \par
{\it We will suppose further that the condition (2) is satisfied:} $ (r,s) \in G_{\alpha,d}. $ \par
 The {\it exact value} of the constant $ K_{\alpha, d}(r,s) = K(r,s), $ i.e. the value

 $$
 V_{\alpha,d}(r,s) = V(r,s) = \sup_{f \in L_r, f \ne 0} \sup_{g \in L_s,g \ne 0}
 \frac{B_{\alpha,d}(f,g)}{|f|_r \ |g|_s } \eqno(3)
 $$
is now known only in the case $ r = s = 2d/(d + \alpha). $  Namely, it is proved
in the articles  \cite{Lieb1}, \cite{Lieb3}  (in our notations) that

$$
V_{\alpha,d}(2d/(d + \alpha), 2d/(d + \alpha)) = \pi^{0.5(d-\alpha)}
\frac{\Gamma(\alpha/2)}{\Gamma(0.5(d+\alpha))} \
\left[\frac{\Gamma(d)}{\Gamma(d/2)} \right]^{\alpha/d}.
$$
 In the book \cite{Lieb0}, p.98 is obtained the following {\it upper} estimate for the
 value $ V_{\alpha,d}(r,s):$

$$
 V_{\alpha,d}(r,s) \le (rs)^{-1} \cdot \alpha^{-1} \cdot
(\omega(d-1))^{1-\alpha/d} \cdot d^{\alpha/d}
\cdot \left[ \frac{r^{1-\alpha/d}}{(r-1)^{1-\alpha/d}} +
\frac{s^{1-\alpha/d}}{(s-1)^{1-\alpha/d}} \right]. \eqno(4)
$$
 Note that the right-hand side of inequality (4) allows a simple estimation:

 $$
 V_{\alpha,d}(r,s)\le C_{1,d} \ \alpha^{-1} \ \left[(r-1)(s-1)\right]^{\alpha/d -1}.
 \eqno(4.a)
 $$

{\sc  We can not calculate the exact values of $ V(r,s)  $ for the different values
$ r,s: r \ne s, $ but we can find the asymptotically equivalent to the right-hand side
as} $ r \to 1+0 $ {\sc and} $ s \to 1 + 0 $
{\sc of the inequality (4.a) low bound  for the value } $ V_{\alpha,d}(r,s) $
{\sc bide-side estimations for  these values  of a view: }

 $$
 V_{\alpha,d}(r,s) \ge
\frac{ K(d) \ \alpha^{-1} }{[(r-1)(s-1)]^{1-\alpha/d }}; \eqno(5)
 $$
here $ (r,s) \in G_{\alpha,d}, \ K(d) $ is finite positive function on the
variable $ d. $ \par
 Evidently, the  estimations (4), (4a) and (5) may be rewritten as follows:

 $$
\frac{ K_3(d) \ \alpha^{-1} }{[(r-1)(s-1)]^{1-\alpha/d }}
\le V_{\alpha,d}(r,s) \le
\frac{ K_4(d) \ \alpha^{-1} }{[(r-1)(s-1)]^{1-\alpha/d }}, \eqno(5.a)
 $$
$ \forall d = 1,2, \ldots \ 0 < K_3(d) \le K_4(d) < \infty.  $\par

{\it  We dare formulate as a hypotheses the following equality: }

$$
\sup_{r,s \in G_{\alpha,d}} \left\{\left[(r-1)^{1-\alpha/d} (s-1)^{1 - \alpha/d} \right]
V_{\alpha,d}(r,s) \right\} = K_5(\alpha,d) \alpha^{-1},
$$
where the powers $ 1-\alpha/d, \ 1 - \alpha/d $ are exact and the function
$ K_5(\alpha,d) $ is positive and continuous in the {\it closed}  interval
 $  \alpha \in [0, d]. $\par

\vspace{4mm}

The article is organized as follows. In the next section we obtain
the main result: low bounds for Riesz's linear operator and correspondent bilinear
form in the Lebesgue spaces. \par
In the third section we consider the case when instead the whole space $ R^d $ in the
integral in (1) is some bounded domain. \par
 The last section contains some slight generalizations of obtained results.\par
\vspace{3mm}

 We use symbols $C(X,Y),$ $C(p,q;\psi),$ etc., to denote positive
constants along with parameters they depend on, or at least
dependence on which is essential in our study. To distinguish
between two different constants depending on the same parameters
we will additionally enumerate them, like $C_1(X,Y)$ and
$C_2(X,Y).$ The relation $ g(\cdot) \asymp h(\cdot), \ p \in (A,B), $
where $ g = g(p), \ h = h(p), \ g,h: (A,B) \to R_+, $
denotes as usually

$$
0< \inf_{p\in (A,B)} h(p)/g(p) \le \sup_{p \in(A,B)}h(p)/g(p)<\infty.
$$
The symbol $ \sim $ will denote usual equivalence in the limit
sense.\par
We will denote as ordinary the indicator function
$$
I(x \in A) = 1, x \in A, \ I(x \in A) = 0, x \notin A;
$$
here $ A $ is a measurable set.\par
 Other notations. We denote as usually

 $$
 \Omega(d) = \frac{\pi^{d/2}}{\Gamma(1 + d/2)}, \ \omega(d) =
 \frac{2 \pi^{d/2}}{\Gamma(d/2)}, \ \overline{\omega}(d)= \max(1,\omega(d)).
 $$
  The value $ \Omega(d) $ is the volume of the unit ball in the space $ R^d: $

  $$
  m(B(x,r)) = \Omega(d) \ r^d; \ B(x,r) = \{y, y \in R^d, |x-y| \le r \},
  $$
 $ m(\cdot) $ denotes ordinary Lebesgue measure and
 $ \Gamma(\cdot) $ denotes usually Gamma - function. \par

  The  value $ \omega(d) $ is the area of the unit sphere in this space. \par
  Note that $ \omega(d) = d \ \Omega(d). $ \par
Further, we need to use the so-called maximal operator $ Mf(x): $
 $$
 M f(x) \stackrel{def}{=} \sup_{r > 0}
 \left[ \Omega(d)^{-1} r^{-d} \ \int_{y: |y - x| \le r} |f(y)| \ dy \right].
 $$
It is known ( \cite{Stein1}, p. 173 - 188) that

$$
| M f|_p \le  C(d) \frac{p}{p-1} \ |f|_p, \ p \in (1, \infty), \ C(d) \in (0,\infty).
$$
The minimal value of the constant $ C(d) $ from the last inequality will be denoted by
$ S(d), $ ("Stein's constant"), on the other hand:

$$
S(d) \stackrel{def}{=} \sup_{p \in (1,\infty)} \sup_{f \ne 0, \ f \in L(1,\infty)}
\frac{|Mf|_p }{p \ |f|_p/(p-1)}.
$$
It is evident that
$$
 \inf_{p \in (1,\infty)} S(d) \ge 1, \ d = 1,2,\ldots.
$$
 The first upper estimation for the value $ S(d) $ was obtained in the classical
 book of E.M.Stein ( \cite{Stein1}, p. 173 - 188):

$$
S(d) \le 2 \cdot 5^{d}.
$$
In the article \cite{Hunt1} it is proved that $ S(2) \le 2. $ In the next works of
E.M.Stein \cite{Stein1}, \cite{Stein2}, \cite{Stein3},\cite{Stein4},it was obtained
the following estimations for $ S(d): $
$$
S(d) \le C_1 \ \sqrt{d}, \ S(d) \le C_2
$$
with {\it some absolute} constants $ C_1 $ and $ C_2. $ \par

 All the passing to the limit in this article may be grounded by means
 of Lebesgue dominated convergence theorem.\par

\vspace{3mm}

{\bf Some comments about upper estimations for Riesz potential.}\par

\vspace{3mm}
 The original proof of the inequality (4) in the book \cite{Lieb0}, pp. 98-112
 based on the {\it weak} Young's inequality

 $$
 \int_{R^d} \int_{R^d} f(x) g(x-y) h(y) \ dx dy \le C(p,q,r) \ |f|_p \ |g|_{q,w} \
 |h|_r,
 $$
where $ p,q,r > 1, 1/p + 1/q + 1/r = 2 $ and $ |g|_{q,w}  $ denotes the weak Lebesgue
norm of order $ q, $ which may be defined  (up to norm equivalence) as
$$
|g|_{q,w} = \sup_{A, m(A) \in (0,\infty)} (m(A))^{1/q-1} \int_A |g(x)| \ dx,
$$
where $ m(A) $ denotes usually Lebesgue measure of (measurable)  set $ A, $
in contradiction to the classical Young's inequality
$$
 \int_{R^d} \int_{R^d} f(x) g(x-y) h(y) \ dx \ dy \le \ |f|_p \ |g|_q \ |h|_r,
$$
where again $ p,q,r > 1, 1/p + 1/q + 1/r = 2. $ \par

{\it But we will offer an another proof, which gives the possibility for many
 generalizations.} \par

\vspace{4mm}

We will consider in this  subsections only the values $ p $ from the open interval
$ p \in (1,d/\alpha) $ and denote $ q = q(p) = pd/(d - \alpha p); $  or equally

$$
\frac{1}{q} = \frac{1}{p} - \frac{\alpha}{d};
$$
evidently, $ q \in (d/(d - \alpha), \infty). $ \par
The inverse function to the function $ q = q(p) $ has a view
$ p = p(q) = dq/(d + \alpha q );  $ note that if $ p \to 1 + 0 \ \Rightarrow
 q \to d/(d - \alpha) + 0 $ and if $ p \to d/\alpha - 0 \ \Rightarrow q \to \infty. $ \par

Note that the case $ p = 1 $ and $ p = d/\alpha, $ i.e. when $ f \in L_1(R^d) $ or
 $ f \in L_{d/\alpha}(R^d) $ is considered, e.g., in \cite{Adams1}, p. 56. \par

 More detail,

 $$
 p - 1 = \frac{(d - \alpha) (q - d/(d - \alpha))}{ d + \alpha \ q}
 $$
and
$$
\frac{d}{\alpha} - p = \frac{d^2}{\alpha(d + \alpha q) }.
$$

{\bf Theorem 1.}
$$
|I_{\alpha}f|_q \le \left[S(d) \ \overline{\omega(d)} \right] \cdot
|f|_p \cdot p^{(\alpha p -d)(p-1)} \cdot(p-1)^{\alpha(p-1)/d } \times
$$

$$
[(p-1)(d/\alpha - p)]^{\alpha/d - 1} \cdot
\left[1 + \frac{(p-1)^{1 - 1/p}}{\alpha p }(d - \alpha p) \right] (\stackrel{def}{=} R), \eqno(6)
$$
which is equivalent the inequality (4). \par
{\bf Remark 1.} Note that the right hand side of the inequality (6), say $ R, $ may be estimated as follows:
$$
R \le \left[S(d) \ \overline{\omega(d)} \right] \cdot
 \ |f|_p \cdot \frac{2 d^2 / \alpha}{[(p-1)(d/\alpha - p)]^{1-\alpha/d    }}.
\eqno(7)
$$
{\bf Remark 2.} The last estimation(7) improved the main result of the article
\cite{Ostrovsky3}.\par

{\bf Proof} of the theorem 1 (briefly).\par
 In the book \cite{Adams1}, pp. 49 - 54 is described a very interest
  approach, which we will use here.\par

Let $ \chi =  \chi(z) $ be a positive at $ z \in (0,\infty) $ continuous decreasing
function such that $ \phi(\infty) \stackrel{def}{=} \lim_{z \to \infty} \phi(z) = 0. $
We define also a function

$$
\Phi(z) = \int_z^{\infty} \chi(t) \ dt,
$$
if there exists. \par
We have  for the values $ \delta \in (0, \infty) $ analogously to the assertion
in \cite{Adams1}, p. 49 - 51:

$$
\int_{y: |x - y| < \delta} \Phi(| x - y |) \ f(y) dy = \int_0^{\delta} \chi(r)
\int_{y: |x - y| < r} f(y) \ dy + \Phi(\delta) \int_{y: |x - y| < \delta} f(y) \ dy.
$$
 Without loss of generality we can assume that the function $ f(\cdot) $ is
 non - negative. \par
 As long as

 $$
 \int_{y: |x - y| \le r} f(y) \ dy \le \Omega(d) \ r^d \ Mf(x),
 $$

 $$
 \int_{y: |x - y| \le \delta} f(y) \ dy \le \Omega(d) \ \delta^d \ Mf(x),
 $$
we obtain the  estimate

$$
\int_{y: |x - y| < \delta} \Phi(|x - y |) \ f(y) dy \le \Omega(d) \ Mf(x) \cdot
\left[\int_0^{\delta} r^d \ \chi(r) \ dr + \delta^d \ \Phi(\delta) \right]
$$

$$
\stackrel{def}{=}  Mf(x) \ A_{\chi, d}(\delta).
$$

Further, we have denoting $ s = p/(p-1) $ and using H\"older's inequality:

$$
\int_{y: |x - y| \ge \delta} \Phi(|x - y|) \ f(y) \ dy \le |f|_p \cdot
\left[ \int_{y: |x - y| > \delta} \Phi^s(|x - y| \ dy  \right]^{1/s} =
$$

$$
\omega(d) \ |f|_p \left[\int_{\delta}^{\infty} r^{d-1} \Phi^s(r) \ dr \right]^{1/s}
= D_{\chi}(p, \delta) \ |f|_p,
$$
where

$$
D_{\chi}(p, \delta) \stackrel{def}{=} \omega(d) \left[\int_{\delta}^{\infty} r^{d - 1} \
\Phi^s(r)  \ dr \right]^{1/s},
$$
if there exists for some values $ s $ from some non - trivial interval
$ s \in (d/(d - \alpha), s_0); $ if $ s_0 < \infty, $ then we define formally
$ D_{\chi}(p, \delta) = + \infty.  $\par

We conclude tacking into account the partition

$$
\int_{R^d} \Phi(|x - y|) \ f(y) \ dy = \int_{y: |x - y| < \delta} \Phi(|x - y|) \ f(y) \ dy
$$

$$
+ \int_{y: |x - y| \ge \delta} \Phi(|x - y|) \ f(y) \ dy:
$$

$$
\int_{R^d} \Phi(|x - y|) \ f(y) \ dy \le Mf(x) \ A_{\chi, d}(\delta) +
D_{\chi}(p, \delta) \ |f|_p.
$$
Therefore, we can make the optimization over $ \delta, \ \delta \in (0,\infty): $

$$
\int_{R^d} \Phi(|x - y|) \ f(y) \ dy \le \inf_{\delta > 0}
\left[ Mf(x) \ A_{\chi, d}(\delta) +D_{\chi}(p, \delta) \ |f|_p \right]
$$

$$
\stackrel{def}{=} H(p, Mf(x), |f|_p).
$$
 Solving the last inequality, we obtain denoting
 $$
 w(x) =  \int_{R^d} \Phi(|x - y|) \ f(y) \ dy:
 $$
the inequality of a view

 $$
 G(p, w(x), |f|_p) \le (Mf(x))^p,
 $$
 and after the integration with the corresponding power

$$
 \int_{R^d} G(p, w(x), |f|_p) \ dx \le |Mf(x)|_p^p \le C^p(\alpha,d) \ |f|_p^p \ (p-1)^{-p},
 $$

 $$
  \left[\int_{R^d} G(p, w(x), |f|_p) \ dx \right]^{1/p}
   \le |Mf(x)|_p \le C(\alpha,d, \chi) \ |f|_p \ (p-1).
 $$

Since the relation between  the functions $ f(\cdot) $ and $ w(\cdot) $ is linear,
the last inequality has a view

$$
|w|_q \le C_2(\alpha,p, \chi(\cdot)) \ |f|_p.
$$

 Choosing the function $ \Phi(r) $ as  follows:

 $$
 \Phi(r)= r^{\alpha - d}
 $$
 and tacking into account the Stein's estimations for the maximal functions $ Mf(x), $
 we obtain the assertion of theorem 1 after  complicate calculations.\par

\vspace{3mm}

\hfill $\Box$

\bigskip

\section{ Main Result: low bounds for Riesz's potential.}

\vspace{3mm}

 In this section we built some examples in order to illustrate the
exactness of upper estimations (4). \par
 We consider here only more hard case $ d \ge 2; $ the one-dimensional case $ d = 1 $
 is considered in the work \cite{Ostrovsky3}. \par
 Before the formulating of the main result of this section, we must introduce
 several new notations. \par

 $$
 a = a(\alpha,d):= e^{1/e}\cdot \max \left[\frac{\omega(d)}{\alpha},
 \left(\frac{\omega(d)}{\alpha} \right)^{d/\alpha}   \right],
 $$

$$
m = m(\alpha,d):= \min \left(1, (\omega(d)/d)^{1-\alpha/d}\right),
$$

$$
A = A(\alpha,d):= 4\pi/(9\alpha) \cdot \omega(d-1)\cdot 2^{-d} \cdot m(\alpha,d),
$$

$$
n = n(\alpha,d):= \max \left(\omega(d)/d, (\omega(d)/d)^{\alpha/d  } \right),
$$

$$
D = D(\alpha,d):= 3^{-1} \cdot 4  \cdot 5^{\alpha-d}\cdot \omega(d) \cdot
\min(1,(\omega(d))^{d/(d-\alpha)} \cdot (d^2/\alpha)^{-2-\alpha/d}.
$$

Let us introduce the following important function: $ F(p) = F_{\alpha,d}(p):= $

$$
0.5 \ \frac{A \cdot(p-1)^{1/p + (d-\alpha)/\alpha} +
D \cdot(d/\alpha - p)^{1/p + (2d-\alpha)/d }}
{a \cdot(p-1)^{1/p} + n \cdot(d/\alpha - p)^{\alpha/d}}. \eqno(8)
$$
 Since the function $ p \to F(p) $ is continuous in the {\it closed} interval
 $ p \in [1,d/\alpha] $ and is positive, we can introduce the strong positive
 variable

 $$
 R(\alpha,d) := \inf_{p \in [1, d/\alpha]} F_{\alpha,d}(p). \eqno(9)
 $$
\vspace{3mm}
{\bf Theorem 2.} For $ (r,s) \in G_{\alpha,d}  $ there holds:
$$
V_{\alpha,d}(r,s) \ge \frac{ R(\alpha,d) }{[(r-1)(s-1)]^{1-\alpha/d}}, \eqno(10)
$$
\vspace{3mm}
{\bf Remark 3.} Note that the $ R(\alpha,d) $ may be estimated from below as follows:
$ R(\alpha,d) \ge K_6(d) \ \alpha^{-1}, \ \alpha \in (0, 1), $  which asymptotically
coincides with the upper bound given by theorem 1.\par
{\bf Proof.}  We will consider two examples of a functions from the set
$ L(1,d/\alpha).$ \par
\vspace{3mm}
{\bf First example.}
$$
f_0(x) = |x|^{-d} \ I(|x| > 1).
$$
We find by direct calculations using multidimensional polar coordinates:

$$
|f_0|_p = (\omega(d))^{1/p} \ d^{-1/p} \ (p-1)^{-1/p} \le n(\alpha,d) \ (p-1)^{-1/p},
p > 1;
$$

$$
u_0(x):= I_{\alpha}f_0(x)\ge C_{\alpha}(d) \ |x|^{\alpha - d} \ |\log |x|| \ I(|x| > 1),
$$
where
$$
C_{\alpha}(d) = 4 \cdot 5^{\alpha - d} \cdot \omega(d);
$$

$$
|u_0|_q \ge C_{\alpha}(d) \cdot (\omega(d))^{1/q} \times
\frac{(\Gamma(q+1))^{1/q}}{(q(d-\alpha)-d)^{1 + 1/q} } \ge
$$

$$
D_{\alpha,d} \cdot (d/\alpha - p)^{1 + 1/q} \cdot (p-1)^{-1-1/q};
$$
recall that $ q = q(p). $ \par
\vspace{3mm}
{\bf Second example.} We put:

$$
g_0(x) = |x|^{-\alpha} \ I(|x| < 1),
$$
and find:

$$
|g_0|_p = \frac{ \omega^{1/p}(d) \ \alpha^{-1/p}}{ (d/\alpha - p)^{1/p}} \le
a(\alpha,d) \cdot (d/\alpha - p)^{-\alpha/d}, \ p \in (1,d/\alpha);
$$

$$
v_0(x) := I_{\alpha}g_0(x) \ge 3^{-1} \cdot 4\pi \cdot \omega(d-1) \cdot  2^{-d}
\cdot |\log |x|| \cdot I(|x| < 1);
$$

$$
|v_0|_q \ge  \omega^{1/q} \cdot d^{-1 - 1/q} \cdot \Gamma^{1/q}(q+1) \ge
A(\alpha,d) \cdot (d/\alpha -p)^{-1}, \ p \in (1, d/\alpha).
$$
\vspace{3mm}
{\bf Third example. Summing.} \par
We define:

$$
h(x) = f_0(x) + g_0(x).
$$
It follows from the triangle inequality:

$$
|h|_p \le |f_0|_p + |g_0|_p \le a(\alpha,d)(d/\alpha - p)^{-\alpha/d} +
n(\alpha,d) (p-1)^{-1/p}.
$$
 Further, as long as both the functions $ u_0 $ and $ v_0 $ are positive, we obtain
 for the values $ p $ from the interval $ p \in (1, d/\alpha): $

 $$
 |I_{\alpha}h|_q \ge 0.5(|u_0|_q + |v_0|_q) \ge A(\alpha,d) (d/\alpha-p)^{-1} +
 D(\alpha,d)(d/\alpha - p)^{1+1/q}(p-1)^{-1-1/q}.
 $$
We get after dividing:

$$
2 \frac{|I_{\alpha}h|_q \cdot [(p-1)(d/\alpha-p)]^{1-\alpha/d}}{ |h|_p  }
\ge \frac{A(p-1)^{1+1/q} + D(d/\alpha-p)^{2 + 1/q}}
{a(p-1)^{1/p} + n(d/\alpha-p)^{\alpha/d}},
$$
i.e.

$$
\frac{|I_{\alpha}h|_q \cdot [(p-1)(d/\alpha-p)]^{1-\alpha/d}}{ |h|_p  } \ge
F_{\alpha,d}(p),
$$
following

$$
\frac{|I_{\alpha}h|_q \cdot [(p-1)(d/\alpha-p)]^{1-\alpha/d}}{ |h|_p  } \ge
R(\alpha,d), \eqno(11)
$$
which is equivalent to the assertion of theorem 2. \par

\bigskip

\section{The case of a bounded domain}

\vspace{3mm}

We consider in this section the  {\it truncated}  Riesz's operator

$$
u^{(G)} = u^{(G)}(x)  = I^{(G)}_{\alpha} f(x) = \int_G \frac{f(x - y) \ dy }{|y|^{d - \alpha} }, \eqno(12)
$$
where $ G $ is open bounded domain in $ R^d $ containing the origin and
 such that

$$
0 < \inf_{x \in \partial G } |x| \le \sup_{x \in \partial G } |x| < \infty,
$$
$ \partial G $ denotes boundary of the set $ G. $ \par
It is known (see, e.g. \cite{Mitrinovich1}, p.90), that if
$ f \cdot \max(1, \log |f|) \in L_1(G), $ then
$ I^{(G)}_{\alpha}f \in L_q(G). $ \par
 We can and will assume further without loss of generality that the set $ G $ is
 unit ball in the space $ R^d: $

 $$
G = \{x, \ x \in R^d, \ |x| \le 1 \}.
 $$
 Let us denote as ordinary $ p^/ = p/(p - 1), p \in (1,\infty) $  and introduce the
 correspondent bilinear truncated Riesz's functional:

$$
B^{(G)}(f,g) = \int_{R^d} g(x) \ dx
 \left \{\int_{G} \frac{f(y) \ dy}{|x - y|^{d - \alpha}} \right\}=
 ( I^{(G)}_{\alpha} f, g). \eqno(14)
 $$
Let us denote
$$
Z(p) = \frac{[\omega(d)/(d-\alpha)]^{1/p}}{(d/(d-\alpha) - p)^{1/p}},
$$
where the variable $ p $ is the following function on the variables $ r $ and $ s: $

$$
p = p(r,s) = \frac{r^/ \ s^/ }{r^/ + s^/ },
$$
and the variable $ r,s $ are such that
$$
r,s \ge 1, \ 1 \le 1/r + 1/s < 1 + \alpha/d.
$$

\vspace{3mm}

{\bf Theorem 3}.
$$
|B^{(G)}(f,g)|\le Z(p(r,s)) |f|_r |g|_s. \eqno(15)
$$
\vspace{3mm}

{\bf Proof } is very simple by means of the {\it classical} Young's inequality
(the {\it weak } Young's inequality gives at the same result). \par

Namely, we obtain by the direct computation for the values $ p $ from the
interval $ p \in [1, d/(d - \alpha) ) $:

$$
| |x|^{\alpha - d} \cdot I(x \in G)|_p = Z(p).
$$

 Note that in contradiction to the unbounded case $ G = R^d $ the values $ (r,s) $
 have two degrees of freedom.\par

\hfill $\Box$

\bigskip

\section{Concluding remarks}

\vspace{3mm}

{\bf A.} We consider in this subsection some generalization of the Riesz's potential (non-truncated) operator of a view

$$
I_{\alpha, \beta} f(x) =
\int_{R^d} \frac{f(y) \  | \log |x - y| \ |^{\beta} \ dy }{ |x - y|^{d - \alpha} },
$$
$ \alpha = \const \in (0,d), \ \beta = \const > 0, $ or equally

$$
I_{\alpha, \beta} f(x) =
\int_{R^d} \frac{f(y) \  [ 1 + | \log |x - y| \ |]^{\beta} \ dy }{ |x - y|^{d - \alpha} },
$$

or more generally

 $$
 I_{\alpha, \beta}^{(Q)} f(x) = \int_{R^d}
 \frac{f(y) \ | \log |x - y| \ |^{\beta} \ Q(|\log|x - y| \ |) \ dy }
 { |x - y|^{d -\alpha} }, \eqno(16)
 $$
where $ \alpha = \const \in (0,d), \ \beta = \const > 0, $ and
$ Q(z) $ is a {\it slowly varying} as $ z \to \infty $ continuous positive function:

$$
\forall \lambda > 0 \ \Rightarrow  \lim_{z \to \infty} Q(\lambda z)/Q(z) = 1.
$$

It is proved in the article \cite{Ostrovsky3} by means of considerations used
in the proof of theorem 1 that

 $$
 |I_{\alpha, \beta}^{(Q)} f|_q \le
 \frac{C \ |f|_p }{[(p - 1) \cdot (d/\alpha - p)]^{1 + \beta - \alpha/d} }
  \cdot Q( q(d - \alpha) - d)^{-1}) \ Q(q), \eqno(17)
 $$
and that the last inequality (17) is asymptotically as $ p \to 1+0 $ and
$ p \to d/\alpha - 0 $ exact, e.g. for the function $ h(\cdot). $ \par
Here as in the sections 1,2 $ p \in (1/d/\alpha) $ and $ q = q(p). $ \par
 Therefore, we have for the corresponding {\it generalized bilinear Riesz's functional}

 $$
 B_{\alpha,d}^{(Q)}(f,g):= (g,I_{\alpha, \beta}^{(Q)} f ),
 $$

 $$
 V_{\alpha,d}^{(Q)}(r,s):= \sup_{f \in L_r, f \ne 0} \sup_{g \in L_s,g \ne 0}
 \frac{B_{\alpha,d}^{(Q)}(f,g)}{|f|_r \ |g|_s }
 $$
the following bilateral estimation: \par
{\bf Theorem 4.}

$$
 \frac{C_1(d,Q(\cdot)) \ \alpha^{-1-\beta + \alpha/d} \ Q(1/\alpha) }
 {[(r - 1) \ (s - 1)]^{1 + \beta - \alpha/d } }
   \le V_{\alpha,d}^{(Q)}(r,s) \le
 \frac{C_2(d,Q(\cdot)) \ \alpha^{-1-\beta + \alpha/d} \ Q(1/\alpha) }
 {[(r - 1) \ (s - 1)]^{1 + \beta - \alpha/d } }, \eqno(18)
 $$
where as in (2)
 $$
 r,s > 1, \ \frac{1}{r} + \frac{1}{s} = 1 + \frac{\alpha}{d}.
 $$

\end{document}